\documentclass{article}

\begin{document}
\textbf{Zermelo's theorem}. Each set can be well ordered.

\textbf{Proof}. I. By $A$ we denote the set under consideration. Let $B$ be the set of all the subsets of $A$.
Let $\phi:B\setminus\{\emptyset\}\to A$ be the function assigning to each nonempty subset $X\subseteq A$ a point
$x\in X$ (by the axiom of choice such a function exists). The function $\alpha(X)=\phi(A\setminus X)$ is defined
for all the subsets of $A$ except $A$ itself.

II. A subset $P\subseteq B$ is said to be \emph{regular} if the following conditions are satisfied:

1) $P$ is linearly ordered with respect to the relation $\subseteq$, i.e. if $p_1,p_2\in P$, then either
$p_1\subseteq p_2$, or $p_2\subseteq p_1$;

2) $P$ is well ordered with respect to the relation $\subseteq$, i.e. if $\gamma\subseteq P$, then $\gamma$ has
a least element in this ordering (note that it is equal to $\cap\gamma$);

3) $\emptyset\in P$;

4) if a set $p \in P$ is not empty, then $p=p_1\cup\{\alpha(p_1)\}$, where $p_1=\cup\{q: q\in P, q\subset p\}$.

Regular sets exist. For example, $\{\emptyset\}$, $\{\emptyset,\{\alpha(\emptyset)\}\}$, $\{\emptyset,
\{\alpha(\emptyset)\}, \{\alpha(\emptyset), \alpha(\{\alpha(\emptyset)\})\}\}$ are such sets. Note that if $p\in
P$ and for the next element we have $p+1\in P$, then $p+1=p\cup\{\alpha(p)\}$.

III. Let $P_1$ and $P_2$ be regular sets. Define

$P_3=\{p: p\in P_1\cap P_2, \{q: q\in P_1, q\subset p\}=\{q: q\in P_2, q\subset p\}\}$.

Let us show that $(*)$ $P_3=P_1$ or $P_3=P_2$.

Suppose the contrary. Since $P_3\subseteq P_1\cap P_2$, it follows that the sets $P_1\setminus P_3$ and
$P_2\setminus P_3$ are not empty. Let $r_1$ be the least element of $P_1\setminus P_3$ and let $r_2$ be the
least element of $P_2\setminus P_3$. Since $\{p: p\in P_1, p\subset r_1\}=P_3=\{p: p\in P_2, p\subset r_2\}$, it
follows from 4) that $r_1=\cup P_3\cup\{\alpha(\cup P_3)\}=r_2$. Hence $r_1\in P_3$, which is impossible
($r_1\notin P_3$).

Thus, having supposed that $(*)$ is not true, we arrive at a contradiction.

So, we have $(*)$, which means that either $P_1$ is an initial segment of $P_2$, or $P_2$ is an initial segment
of $P_1$.

IV. Let us denote by $Q$ the union of all the regular sets. The set $Q$ obviously satisfies the conditions 1)
and 3) from II.

Let us show that 2) holds.

Let $\gamma\subseteq Q$. For some regular set $P$ the intersection $\gamma\cap P$ is not empty. Let
$m\in\gamma\cap P$. By virtue of regularity of $P$ the set $\{n: n\in\gamma, n\subseteq m\}\subseteq P$ has a
least element. We denote it by $g$. By the definition it is not greater than any element of $\gamma$ less than
$m$, and it is not greater than any element of $\gamma$ greater than $m$, since $g\subseteq m$.

Let us show that 4) holds.

Suppose that $q\in Q$ is not empty. By the definition of $Q$ there is a regular set $P$ such that $q\in P$. By
4) we have $q=q_1\cup\{\alpha(q_1)\}$ where $q_1=\cup\{p: p\in P, p\subset q\}$. For any set $r\in Q\setminus P$
we have $q\subseteq r$, so $q_1=\cup\{p: p\in Q, p\subset q\}$.

Thus, $Q$ is a regular set. Let $Z=\cup Q$.

If $Z\neq A$, then the set $\tilde{Q}=Q\cup\{Z\cup\{\alpha(Z)\}\}$ is regular. It contradicts the definition of
$Q$ as the union of all the regular sets, since $\tilde{Q}$ contains $Q$ as a proper subset.

Thus, $\cup Q=A$.

V. Consider $\alpha$ as a map from $Q$ to $A$.

Let us show that $\alpha$ is injective.

Let $q_1\neq q_2$. Without loss of generality we may assume that $q_1 \subset q_2$. Then $q_1+1\subseteq q_2$.
Since $Q$ is regular, we have $q_1+1=q_1\cup\{\alpha(q_1)\}$. Therefore, $\alpha(q_1)\in q_1+1\subseteq q_2$,
i.e. $\alpha(q_1)\in q_2$. But $\alpha(q_2)\notin q_2$. Hence $\alpha(q_1)\neq\alpha(q_2)$.

Let us show that $\alpha$ is surjective.

As $\cup Q=A$ for every $a\in A$, the set $M_a=\{q: q\in Q, q\ni a\}$ is not empty. Denote by $r$ the least
element of $M_a$. By regularity of $Q$ we have $r=r_1\cup\{\alpha(r_1)\}$, where $r_1=\cup\{q: q\in Q, q\subset
r\}$. Since $r$ is the least element containing $a$, we have $a\notin r_1$. Hence $\alpha(r_1)=a$.

Thus, $\alpha$ induces a well-order relation on $A$. The proof is completed.
\end{document}